\documentstyle[12pt, amstex]{article}

in \oddsidemargin 0true in \evensidemargin 0true in \headsep
0.3true in

\begin{document}
\newcommand{\bx}{\hfill\rule{.25cm}{.25cm}\medbreak}
\newtheorem{prop}{Proposition}[section]
\newtheorem{defn}{Definition}[section]
\newtheorem{lem}{Lemma}
\newtheorem{thm}{Theorem}[section]
\newtheorem{asum}{Assumption}
\newtheorem{cor}{Corollary}

\newcommand{\bR}{\overline{R}}
\newcommand{\eg}{\frak g}
\newcommand{\stu}{{\frak stu}_n(R, -, \boldsymbol{\gamma})}
\newcommand{\eij}{e_{ij}(a)}
\newcommand{\ejk}{e_{jk}(b)}
\newcommand{\Tij}{T_{ij}(a, b)}
\newcommand{\eik}{e_{ik}(ab)}
\newcommand{\ekl}{e_{kl}(b)}
\newcommand{\st}  {{\frak st}(m,n,R)}
\newcommand{\sth}  {\widehat{{\frak st}}_n(R)}
\newcommand{\stf} {{\frak st}(3,1,R)}
\newcommand{\stfh} {\widehat{{\frak st}}(3,1,R)}
\newcommand{\stft} {\widetilde{{\frak st}}(3,1,R)}
\newcommand{\stfs} {{\frak st}(3,1,R)^\sharp}
\newcommand{\stg} {{\frak st}(2,2,R)}
\newcommand{\stgh} {\widehat{{\frak st}}(2,2,R)}
\newcommand{\stgt} {\widetilde{{\frak st}}(2,2,R)}
\newcommand{\stgs} {{\frak st}(2,2,R)^\sharp}
\newcommand{\stt} {{\frak st}(2,1,R)}
\newcommand{\stth} {\widehat{{\frak st}}(2,1,R)}
\newcommand{\sttt} {\widetilde{{\frak st}}(2,1,R)}
\newcommand{\stts} {{\frak st}_3(R)^\sharp}
\newcommand{\lij}{\gamma_i\gamma_j^{-1}}
\newcommand{\lji}{\gamma_j\gamma_i^{-1}}
\newcommand{\ljk}{\gamma_j\gamma_k^{-1}}
\newcommand{\gd}{\dot{\frak{g}}}

\begin{center}{\large \bf Central extensions of Steinberg Lie
superalgebras of small rank}\end{center}

\vspace{1em} \centerline{\bf Hongjia Chen, Yun Gao and Shikui
Shang\footnote[1] {The corresponding author.\\
 Research of the second
author was partially supported by NSERC of
 Canada and Chinese Academy of Science. } }

\vspace{1.5em}
\begin{abstract}
It was shown   by A.V.Mikhalev and I.A.Pinchuk in [MP] that the
second homology group $H_2(\st)$ of the Steinberg Lie superalgebra
$\st$ is trivial for $m+n\geq 5$. In this paper, we will work out
$H_2(\st)$ explicitly for $m+n=3, 4$.
\end{abstract}

\vspace{0.5em}

\noindent{\bf Introduction}

 Steinberg Lie algebras ${\frak st}_n(R)$ play an important role in
 (additive) algebraic K-theory. They have been studied by many
 people (see [L] and [GS], and the references therein). The point is
 that for any unital associative algebra $R$ over a field the
 Steinberg Lie algebra ${\frak st}_n(R)$ is the universal central
 extension of $sl_n(R)$ with the kernel isomorphic to the first
 cyclic homology group $HC_1(R)$ except when both $n$ and the characteristic of the field are
 small. As seen in [GS],
if $n=3,4$, $H_2({\frak st}_n(R))$ is not necessarily equal to $0$.

Recently, A.V.Mikhalev and I.A.Pinchuk [MP] studied the Steinberg
Lie superalgebras $\st$ which are central extensions of Lie
superalgebras $sl(m, n, R)$. They further showed that when $m+n\geq
5$,  $\st$ is the universal central extension of $sl(m,n,R)$ whose
kernel is isomorphic to $(HC_1(R))_{\bar{0}}\oplus(0)_{\bar{1}}$,
here we would like to emphasize the ${\mathbb Z}_2$-gradation of the
kernel.

In this paper,  we shall work out $H_2(\st)$ explicitly for $m+n=3,
4$ without any assumption on char $K$ by adopting the definition for
Lie superalgebras (including char $K=2$ case) introduced by Neher
[N]. It is equivalent to work on the Steinberg Lie superalgebras
$\st$ for small $m+n$. This completes the determination of the
universal central extensions of the Lie superalgebras $\st$ and
$sl(m,n,R)$ as well.

For any non-negative integer $m$, set
$${\mathcal I}_m = mR + R[R, R] \ \ \  \text{ and } \ \ \ R_m =
R/{\mathcal I}_m.$$

Our main result of this paper is the following.

\newtheorem{mthm}{Main Theorem}
\renewcommand{\themthm}{}
\begin{mthm} Let $K$ be a unital commutative ring and $R$ be a
unital associative $K$-algebra. Assume that $R$ has a $K$-basis
containing the identity element. Then
$$H_2(\stt)= (0);$$
$$H_2(\stf)= (0)_{\bar{0}}\oplus (R_2^6)_{\bar{1}};$$
$$H_2(\stg)= (R_2^4\oplus R_0^2)_{\bar{0}}\oplus(0)_{\bar{1}}$$
where $R_m^q$ is the direct sum of $q$ copies of $R_m$.
\end{mthm}

It may be noteworthy to point out that $H_2(\stt)= (0)$ unlike the
Lie algebra case in which $H_2(\frak{st}_3(R))$ is not necessarily
zero.

The organization of this paper is as follows. In Section 1, we
review some basic facts on Steinberg Lie superalgebras $\st$.
Section 2 will treat the $m=2,n=1$ case. Section 3 and 4 will handle
the $m=3,n=1$ case and the $m=2,n=2$ case respectively. Finally in
Section 5 we make a few concluding remarks.

\vspace{4mm}

\noindent{\bf \S 1 Basics on $\st$}

Let $K$ be a unital commutative ring. The following definition was
given in [N].

 \noindent{\bf Definition} A $K$-superalgebra
$L=L_{\bar{0}}\oplus L_{\bar{1}}$ with product $[ , ]$ is a Lie
superalgebra if for any homogenous  $x, y, z\in L, w\in
L_{\bar{0}}$,
\begin{align}
 &[y,x]=-(-1)^{{\text deg}(x){\text deg}(y)}[y,x] \tag{S1}\\
 &[x,[y,z]]=[[x,y],z]+(-1)^{{\text deg}(x){\text deg}(y)}[y,[x,z]]\tag{S2}\\
 &[w,w]=0.  \tag{S3}
\end{align}

Note that (S3) is not needed if char $K\neq 2$ and this definition
works well for $K$ of any characteristic.

Let $R$ be a unital associative $K$-algebra. We always assume that
$R$ has a $K$-basis $\{r_\lambda\}_{\lambda\in\Lambda}$($\Lambda$ is
an index set), which contains the identity element $1$ of $R$, i.e.
$1\in\{r_\lambda\}_{\lambda\in\Lambda}$.

 $\Omega=\{1,\ldots,m,m+1,\ldots,m+n\}$ has a
partition $\Omega=\Omega_0\biguplus\Omega_1$, where
$\Omega_0=\{1,\ldots,m\}$ and $\Omega_1=\{m+1,\ldots,m+n\}$. We
define a map $\omega:\Omega\rightarrow{\mathbb Z}_2$, such that
$$\omega(i)=\cases \bar 0 &\text{ for }i\in \Omega_0\\
\bar 1 &\text{ for }i\in\Omega_1\endcases$$

  The $K$-Lie superalgebra of $(m+n)\times (m+n)$ matrices with coefficients in $R$ is denoted by
$gl(m,n,R)$, such that $deg(e_{ij}(a))=\omega(i)+\omega(j)$ for
$a\in R$, $1\leq i, j\leq m+n$. For $m+n\geq 3$, the elementary Lie
superalgebra $sl(m,n,R)$ is the subalgebra of $gl(m,n,R)$ generated
by the elements $e_{ij}(a)$, $1\leq i\neq j\leq m+n$. Note that
$sl(m,n,R)$ can be equivalently defined as $sl(m,n,R) =[gl(m,n,R),
gl(m,n,R)]$, the derived subalgebra of $gl(m,n,R)$, or
$sl(m,n,R)=\{X\in gl(m,n,R)| \text{ str}(X) \in [R, R]\}$, where
$\text{str}(X)$ is the supertrace of $X=(x_{ij})\in M_{m+n}(R)$
given by
$\text{str}(X)=\sum_{i=1}^{m}x_{ii}-\sum_{j=m+1}^{m+n}x_{jj}$.

Clearly, for any $a, b\in R$,
\begin{equation} [\eij , \ejk ]
= \eik  \notag \tag{1.1}\end{equation} if $i, j, k$ are distinct
and
 \begin{equation} [\eij , \ekl ] = 0 \tag{1.2}\end{equation}
if $j\neq k, i \neq l$.

For $m+n\geq 3$, the Steinberg Lie superalgebra $\st$ is defined to
be the Lie superalgebra over $K$ generated by the homogeneous
elements $X_{ij}(a)$, with deg$(X_{ij}(a))=\omega(i)+\omega(j)$ for
any $a\in R$, $1\leq i\neq j\leq m+n$, subject to the relations(see
[MP]):
\begin{align} &a\mapsto X_{ij}(a) \text{ is a $K$-linear map,}\tag{1.3}\\
&[X_{ij}(a), X_{jk}(b)] = X_{ik}(ab), \text{ for distinct } i, j, k, \tag{1.4}\\
&[X_{ij}(a), X_{kl}(b)] = 0, \text{ for } j\neq k, i\neq l,
\tag{1.5}
\end{align}
where $ a, b\in R,$ $ 1\leq i, j, k, l \leq m+n$.

Both Lie superalgebras $sl(m,n,R)$ and $\st$ are perfect(a Lie
superalgebra $\eg$ over $K$ is called perfect if $[\eg , \eg ] =
\eg $). The Lie superalgebra epimorphism:
\begin{equation} \phi: \st \to sl(m,n,R), \tag{1.6}\end{equation}
such that $\phi(X_{ij}(a)) = \eij$, is a central extension and the
kernel of $\phi$ is isomorphic to $HC_1(R)$ (called $HC_2(R)$ in
[MP]), which is the first cyclic homology group of $R$. Eventually,
$HC_1(R)$ is the even part of ker$(\phi)$ and the odd part is equal
to $0$. So the universal central extension of $sl(m,n,R)$ is also
the universal central extension of $\st$ denoted by
$\widehat{\frak{st}}(m,n,R)$. Our purpose is to calculate
$\widehat{\frak{st}}(m,n,R)$ for any ring $K$ and $m+n\geq 3$.

Setting
\begin{equation} \Tij = [ X_{ij}(a), X_{ji}(b)],
\tag{1.7}\end{equation}
\begin{equation} t(a, b) = T_{1j}(a, b) - T_{1j}(1, ba),\tag{1.8}
\end{equation} for $a, b\in R, 1\leq i\neq j\leq m+n$. Both $\Tij$ and $t(a,b)$
are even elements. Then $t(a,b)$ does not depend on the choices of
$j$(see [MP]). Note that $\Tij$ is $K$-bilinear, and so is
$t(a,b)$.

\newtheorem{newlemma}{Lemma 1.9}
\renewcommand{\thenewlemma}{}
\begin{newlemma}For any $a,b,c\in R$, and distinct $i,j,k$, we
have
\begin{align} & T_{ij}(a, b)= -(-1)^{\omega(i)+\omega(j)}T_{ji}(b, a)\tag{1.10} \\
&[T_{ij}(a, b), X_{kl}(c)]=0 \text{ for distinct }i, j, k, l\tag{1.11} \\
&[T_{ij}(a, b), X_{ik}(c)]=X_{ik}(abc),\ \ [T_{ij}(a, b), X_{ki}(c)]= -X_{ki}(cab) \tag{1.12} \\
&[T_{ij}(a, b), X_{jk}(c)]=-(-1)^{\omega(i)+\omega(j)}X_{jk}(bac),\
 \ [T_{ij}(a, b), X_{kj}(c)]=(-1)^{\omega(i)+\omega(j)}X_{kj}(cba) \tag{1.13} \\
&[T_{ij}(a, b), X_{ij}(c)]= X_{ij}(abc + (-1)^{\omega(i)+\omega(j)}cba) \tag{1.14}  \\
&[t(a, b), X_{1i}(c)]= X_{1i}((ab-ba)c), \ \ [t(a, b), X_{i1}(c)] = -X_{i1}(c(ab-ba)) \tag{1.15} \\
 &[t(a, b), X_{jk}(c)]=0
\text{ for }j, k\geq 2\tag{1.16}
\end{align}
\end{newlemma}
{\bf Proof:} By super-antisymmetry, one has:
$$T_{ij}(a, b)=-(-1)^{\omega(i)+\omega(j)}[X_{ji}(b),X_{ij}(a)]=-(-1)^{\omega(i)+\omega(j)}T_{ji}(b,a)$$
From the super-Jacobi identity, we have
$$[x,[y,z]]=[[x,y],z]+(-1)^{{\text deg}(x){\text
deg}(y)}[y,[x,z]]$$ which is equivalent to
$$[[x,y],z]=[x,[y,z]]+(-1)^{{\text deg}(y){\text
deg}(z)}[[x,z],y].$$

So (1.11) is obvious, and
\begin{align*}[T_{ij}(a, b), X_{ik}(c)]&=[[X_{ij}(a),X_{ji}(b)],X_{ik}(c)]=[X_{ij}(a),[X_{ji}(b),X_{ik}(c)]]=X_{ik}(abc)\\
[T_{ij}(a, b), X_{ki}(c)]&=[[X_{ij}(a),X_{ji}(b)],X_{ki}(c)]=(-1)^{(\omega(i)+\omega(j))(\omega(k)+\omega(i))}[[X_{ij}(a),X_{ki}(c)],X_{ji}(b)]\\
&=-(-1)^{2(\omega(i)+\omega(j))(\omega(k)+\omega(i))}[[X_{ki}(c),X_{ij}(a)],X_{ji}(b)]=-X_{ki}(cab)
\end{align*}
which gives (1.12).

Replaced $T_{ij}(a, b)$ by
$-(-1)^{\omega(i)+\omega(j)}T_{ji}(b,a)$ and exchanging $i$ and
$j$, we can obtain (1.13) from (1.12).

For (1.14), we have
\begin{align*}
[T_{ij}(a,b),X_{ij}(c)]&=[T_{ij}(a,b),[X_{ik}(c),X_{kj}(1)]]\\
&=[[T_{ij}(a, b),X_{ik}(c)],X_{kj}(1)]+[X_{ik}(c),[T_{ij}(a,b),X_{kj}(1)]]\\
&=[X_{ik}(abc),X_{kj}(1)]+(-1)^{\omega(i)+\omega(j)}[X_{ik}(c),X_{kj}(ba)]\\
&=X_{ij}(abc+(-1)^{\omega(i)+\omega(j)}cba).
\end{align*}

From (1.8) we obtain
\begin{align*}
[t(a, b),X_{1i}(c)]&=[T_{1j}(a,
b),X_{1i}(c)]-[T_{1j}(1,ba),X_{1i}(c)]\\
&=X_{1i}(abc)-X_{1i}(bac)=X_{1i}((ab-ba)c)
\end{align*}
and $[t(a, b), X_{i1}(c)] = -X_{i1}(c(ab-ba))$, which show that
(1.15) holds true.

(1.16) is easy and the proof is completed.  $\square$

By the above Lemma, we have

\newtheorem{fprop}{Lemma 1.17}
\renewcommand{\thefprop}{}
\begin{fprop} Let ${\frak T}: = \sum_{1\leq i < j\leq m+n}[ X_{ij}(R), X_{ji}(R)]$.
Then ${\frak T}$ is a subalgebra of $\st$ containing the center
$\frak Z$ of $\st$ with $[{\frak T}, X_{ij}(R)] \subseteq
X_{ij}(R)$. Moreover,
\begin{equation}\st = {\frak T}\oplus_{1\leq i \neq j\leq m+n}X_{ij}(R).\tag{1.18}
\end{equation}
\end{fprop}

As for the decomposition of $\st$, we take
$\{r_\lambda\}_{\lambda\in\Lambda}$, the fixed $K$-basis of $R$,
then $\{X_{ij}(r)\}$
$({r\in\{r_\lambda\}_{\lambda\in\Lambda},1\leq i\neq j\leq m+n})$
can be extended to a $K$-basis $\Gamma$ of $\st$.

In fact, the subalgebra $\frak T$ has a more refined structure.

One  can  easily prove the following lemma (see [MP]).
\newtheorem{blem}{Lemma 1.19}
\renewcommand{\theblem}{}
\begin{blem}
Every element $x\in \frak{T}$ can be written as
$$ x = \sum_{i}t(a_i, b_i) + \sum_{2\leq j\leq m+n} T_{1j}(1, c_j),$$
where $a_i, b_i, c_j\in R$.
\end{blem}

The following result is known(see [MP, Theorem 2]).
\newtheorem{blthm}{Theorem 1.20}
\renewcommand{\theblthm}{}
\begin{blthm} {\normalshape} If  $m+n\geq 5$,  then $\phi: \st \to sl(m,n,R)$ gives the universal
central extension of $sl(m,n,R)$ and  the second homology group of
Lie superalgebra $\st$ is $H_2(\st)= 0$.
\end{blthm}

\vspace{0.5em}

\vspace{4mm}

\noindent{\bf \S 2 Central extensions of $\stt$}

In this section we shall treat $H_2(\stt)$.
\newtheorem{tthm}{Theorem 2.1}
\renewcommand{\thetthm}{}
\begin{tthm} {\normalshape} $H_2(\stt)=0$,
i.e. $\stt$ is centrally closed.
\end{tthm}

\noindent{\bf Proof:} Suppose that
\begin{equation}0\rightarrow{\mathcal
V}\rightarrow \sttt\overset{\tau}{\rightarrow}\stt\rightarrow\notag
0\end{equation} is a central extension of $\stt$. We must show that
there exists a Lie superalgebra homomorphism
$\eta:\stt\rightarrow\sttt$ so that $\tau\circ\eta={\text id}$.

Using the $K$-basis $\{r_\lambda\}_{\lambda\in\Lambda}$ of $R$, we
choose a preimage $\widetilde{X}_{ij}(a)$ of $X_{ij}(a)$ under
$\tau$, $1\leq i\neq j\leq 3,a\in
\{r_\lambda\}_{\lambda\in\Lambda}$.
 Let
$\widetilde{T}_{ij}(a,b)=[\widetilde{X}_{ij}(a),\widetilde{X}_{ji}(b)]$,
then
$$[\widetilde{T}_{ik}(1,1),\widetilde{X}_{ij}(a)]=\widetilde{X}_{ij}(a)+{\mu}_{ij}(a)$$
where ${\mu}_{ij}(a)\in{\mathcal V}$ and $i,j,k$ are distinct.
Replacing $\widetilde{X}_{ij}(a)$ by
$\widetilde{X}_{ij}(a)+\mu_{ij}(a)$, then the elements
$\widetilde{X}_{ij}(b)$ still satisfy the relations (1.3). By
super-Jacobi identity (S2), we have
$$\left[\widetilde{T}_{ik}(1,1),[\widetilde{X}_{ik}(a),\widetilde{X}_{kj}(b)]\right]=\left[[\widetilde{T}_{ik}(1,1),\widetilde{X}_{ik}(a)],\widetilde{X}_{kj}(b)\right]+
\left[\widetilde{X}_{ik}(a),[\widetilde{T}_{ik}(1,1),\widetilde{X}_{kj}(b)]\right]
$$ which yields
\begin{align*}
[\widetilde{T}_{ik}(1,1),\widetilde{X}_{ij}(ab)]&=[\widetilde{X}_{ik}(a+(-1)^{\omega(i)+\omega(k)}a),\widetilde{X}_{kj}(b)]-(-1)^{\omega(i)+\omega(k)}[\widetilde{X}_{ik}(a),\widetilde{X}_{kj}(b)]\\
&=[\widetilde{X}_{ik}(a),\widetilde{X}_{kj}(b)].
\end{align*}

We thus have
\begin{equation}
\widetilde{X}_{ij}(ab)=[\widetilde{X}_{ik}(a),\widetilde{X}_{kj}(b)].\tag{2.2}
\end{equation}
For $k\neq i , k\neq j$, we have
\begin{align}
&[\widetilde{X}_{ij}(a),\widetilde{X}_{ij}(b)]\\
=&\left[\widetilde{X}_{ij}(a),[\widetilde{X}_{ik}(b),\widetilde{X}_{kj}(1)]\right]\notag\\
=&\left[[\widetilde{X}_{ij}(a),\widetilde{X}_{ik}(b)],\widetilde{X}_{kj}(1)\right]+
\left[\widetilde{X}_{ik}(b),(-1)^{(\omega(i)+\omega(j))(\omega(i)+\omega(k))}[\widetilde{X}_{ij}(a),\widetilde{X}_{kj}(1)]\right]\notag\\
=&0+0=0.\tag{2.3}
\end{align}

Next, we show that both of
$[\widetilde{X}_{ij}(a),\widetilde{X}_{ik}(b)]$ and
$[\widetilde{X}_{ij}(a),\widetilde{X}_{kj}(b)]$ are equal to $0$.

Since there  is always one element between $\widetilde{X}_{ij}(a)$
and $\widetilde{X}_{ik}(b)$ which is odd, we can assume that it is
$\widetilde{X}_{ij}(a)$, i.e. $\omega(i)+\omega(j)=\bar{1}$ , then
\begin{align}
0&=\left[\widetilde{T}_{ij}(1,1),[\widetilde{X}_{ij}(a),\widetilde{X}_{ik}(b)]\right]\\
&=\left[[\widetilde{T}_{ij}(1,1),\widetilde{X}_{ij}(a)],\widetilde{X}_{ik}(b)\right]
+\left[\widetilde{X}_{ij}(a),[\widetilde{T}_{ij}(1,1),\widetilde{X}_{ik}(b)]\right]\notag\\
&=0+[\widetilde{X}_{ij}(a),\widetilde{X}_{ik}(b)]=[\widetilde{X}_{ij}(a),\widetilde{X}_{ik}(b)]\tag{2.4}
\end{align}
The other cases are similar. Therefore we have
\begin{align}[X_{ij}(a),
X_{kl}(b)] = 0, \text{ for } j\neq k, i\neq l,  a, b\in R, 1\leq i,
j, k, l \leq 3.\tag{2.5}
\end{align}

By our choices, we know that $\widetilde{X}_{ij}(a)$ satisfy the
relation (1.3)-(1.5). Since we have (2.2) and (2.5), by universal
property of $\stt$ there exists a (unique) Lie superalgebra
homomorphism
$$\eta:\stt\rightarrow\sttt$$
such that $\eta(X_{ij}(a))=\widetilde{X}_{ij}(a)$. Evidently,
 $\tau\circ\eta={\text id}$ which implies that the original
 sequence splits. So $\stt$ is centrally closed.
 $\Box$

\noindent {\bf Remark 2.6}  This result is very different from the
one of ${\frak st}_3(R)$(See [GS]). In that case, $H_2({\frak
st}_3(R))=R_3^6$.

\vspace{4mm}

\noindent{\bf \S 3 Central extensions of $\stf$}

In this section, we shall compute the universal central extension
$\stfh$ of $\stf$.

We don't put any assumption on the characteristic of $K$.

For any nonnegative integer $m$, let ${\mathcal I}_m$ be the ideal
of $R$ generated by the elements: $ma$ and $ab-ba$, for $a,b \in R$.
Immediately, we have ([GS, Lemma 2.1])
\newtheorem{flem}{Lemma 3.1}
\renewcommand{\theflem}{}
\begin{flem}
${\mathcal I}_m=mR+R[R,R] \text{ and } [R, R]R=[R, R]R.$
\end{flem}

Let $$R_m:=R/{\mathcal I}_m$$ be the quotient algebra over $K$
which is commutative.  Write $\bar a=a+{\mathcal I}_m$ for $a\in
R$. Note that if $m=2$,  $\overline{a}=-\overline{a}$ in $R_m$.

\newtheorem{fdef}{Definition 3.2}
\renewcommand{\thefdef}{}
\begin{fdef}
${\mathcal W}=R_2^6$ is the direct sum of six copies of $R_2$ and
$\epsilon_{m}(\overline{a})=(0,\cdots,\overline{a},\cdots,0)$ is
the element of $\mathcal W$, of which the $m$-th component is
$\overline{a}$ and others are zero, for $1\leq m\leq 6$.
\end{fdef}

Let $S_4$ be the symmetric group of $\{1,2,3,4\}$.
$$P=\{(i,j,k,l)|\{i,j,k,l\}=\{1,2,3,4\}\}$$ is the set of all the
quadruple with the  distinct components. $S_4$ has a natural
transitive action on $P$ given by
$\sigma((i,j,k,l))=(\sigma(i),\sigma(j),\sigma(k),\sigma(l))$, for
any $\sigma\in S_4$. $$H=\{(1),(13),(24),(13)(24)\}$$ is a
subgroup of $S_4$ with $[S_4:H]=6$. Then $S_4$ has a partition of
cosets with respect to $H$, denoted by
$S_4=\bigsqcup_{m=1}^6\sigma_mH$. We can obtain a partition of
$P$, $P=\bigsqcup_{m=1}^6P_m$, where $P_m=(\sigma_mH)((1,2,3,4))$.
We define the index map
$$\theta:P\rightarrow\{1,2,3,4,5,6\}$$ by
$$\theta\left((i,j,k,l)\right)=m \text{ if } (i,j,k,l)\in P_m,$$ for $1\leq
m\leq 6$.

Using the decomposition (1.18) of $\stf$, we take a $K$-basis
$\Gamma$ of $\stf$, which contains
$\{X_{ij}(r)|r\in\{r_\lambda\}_{\lambda\in\Lambda},1\leq i\neq
j\leq4\}$. Define $\psi: \Gamma \times \Gamma \to \mathcal W$ by
$$\psi(X_{ij}(r),X_{kl}(s))=\epsilon_{\theta((i,j,k,l))}(\overline{rs})\in\mathcal
W,$$ for $r,s\in\{r_\lambda\}_{\lambda\in\Lambda}$ and distinct
$i, j, k, l$ and $\psi=0$, otherwise. Then we obtain the
$K$-bilinear map $\psi:\stf\times\stf\to\mathcal W$ by linearity.

Recall that a Lie superalgebra over $K$ is defined to be an
${\mathbb Z}_2$-graded algebra satisfying $[x, y]=-(-1)^{{\text
deg}(x){\text deg}(y)}[y,x]$,
$$(-1)^{{\text deg}(x){\text
deg}(z)}[[x,y], z] + (-1)^{{\text deg}(x){\text deg}(y)}[[y, z], x]
+ (-1)^{{\text deg}(y){\text deg}(z)}[[z, x], y]=0$$ and $[w,w]=0$
for the homogenous elements $x,y,z\in L$ and $w\in L_{\bar{0}}$.

 We now have
\newtheorem{slemma}{Lemma 3.3}
\renewcommand{\theslemma}{}
\begin{slemma} The bilinear map $\psi$ is a (super) $2$-cocycle.
\end{slemma}
{\bf Proof: }A bilinear map $\psi$ is called a (super) $2$-cocycle,
if it is (super) skew-symmetric and
$$(-1)^{{\text deg}(x){\text
deg}(z)}\psi([x,y], z) + (-1)^{{\text deg}(x){\text deg}(y)}\psi([y,
z], x) + (-1)^{{\text deg}(y){\text deg}(z)}\psi([z, x], y)=0$$ for
homogenous elements $x,y, z\in L$ and $\psi(w,w)=0$ for $w\in
L_{\bar{0}}$.

Since $R_2=R/{\mathcal I}_2$,
$\overline{ab}=\overline{a}\overline{b}=\overline{b}
\overline{a}=\overline{ba}$ and $\overline{a}=-\overline{a}$  for
$a,b\in R$. Thus the order of factors and $\pm$ sign don't play any
role. We can follow the same arguments as in [GS, Lemma 2.3] for
Steinberg Lie algebra $st_4(R)$ to complete the proof. $\Box$

Since $${\mathcal W}=span_{K}\{\psi(X_{ij}(a),X_{kl}(b))|a,b\in R
{\text \ and \ } i,j,k,l {\text \ are \ distinct}\}$$ and
$$\omega(i) + \omega(j) + \omega(k) + \omega(l) = \bar{1} \text{ for distinct } i, j, k, l,$$
 we obtain a central extension of Lie superalgebra
$\stf$, satisfying that $\mathcal W$ is the odd part of the kernel
:
\begin{equation}
0\rightarrow(0)_{\bar{0}}\oplus({\mathcal
W})_{\bar{1}}\rightarrow\stfh\overset{\pi}\rightarrow\stf\rightarrow
0,\tag{3.4}
\end{equation}
i.e. \begin{equation} \stfh=\left((0)_{\bar{0}}\oplus({\mathcal
W})_{\bar{1}}\right)\oplus\stf,\tag{3.5}
\end{equation}with bracket
$$[(c,x),(c',y)]=(\psi(x,y),[x,y])$$
for all $x,y\in\stf$ and $c,c'\in{\mathcal W}$, where $\pi:{\mathcal
W}\oplus\stf\rightarrow\stf$ is the second coordinate projection
map. Then, $(\stfh,\pi)$ is  a central extension of $\stf$. We will
show that $(\stfh,\pi)$ is the universal central extension of
$\stf$. To do this, we define a Lie superalgebra $\stfs$ to be the
Lie superalgebra generated by the symbols $X_{ij}^{\sharp}(a)$,
$i\neq j,a\in R$ and the $K$-linear space ${\mathcal W}$, with
deg$(X_{ij}^{\sharp}(a))=\omega(i)+\omega(j)$ and deg$(w)=\bar{1}$
for any $w\in{\mathcal W}$, satisfying the following relations:
\begin{align} &a\mapsto X_{ij}{^\sharp}(a) \text{ is a $K$-linear mapping,}\tag{3.6}\\
&[X_{ij}^{\sharp}(a), X_{jk}^{\sharp}(b)] = X_{ik}^{\sharp}(ab), \text{ for distinct } i, j, k, \tag{3.7}\\
&[X_{ij}^{\sharp}(a),{\mathcal W}]=0, \text{ for distinct } i, j, \tag{3.8}\\
&[X_{ij}^{\sharp}(a),X_{ij}^{\sharp}(b)]=0, \text{ for distinct } i, j, \tag{3.9}\\
&[X_{ij}^{\sharp}(a),X_{ik}^{\sharp}(b)]=0, \text{ for distinct } i, j, k, \tag{3.10}\\
&[X_{ij}^{\sharp}(a),X_{kj}^{\sharp}(b)]=0, \text{ for distinct } i, j, k, \tag{3.11}\\
&[X_{ij}^{\sharp}(a),
X_{kl}^{\sharp}(b)]=\epsilon_{\theta((i,j,k,l))}(\overline{ab}),
\text{ for distinct } j, k, i, l, \tag{3.12}
\end{align}
where $a,b\in R$,$1\leq i,j,k, l\leq 4$. As $1\in R$, $\stfs$ is
perfect. Clearly, there is a unique Lie superalgebra homomorphism
$\rho:\stfs\rightarrow\stfh$ such that
$\rho(X^\sharp_{ij}(a))=X_{ij}(a)$ and $\rho|_{\mathcal W}=id$.

\noindent{\bf Remark 3.13: }Comparing with the relations of $\st$
(1.3)-(1.5), we separate the case $[X_{ij}^{\sharp}(a),
X_{kl}^{\sharp}(b)] ( j\neq k, i\neq l)$ into four subcases
(3.9)-(3.12).

We claim that $\rho$ is actually an isomorphism.

\newtheorem{tlemma}{Lemma 3.14}
\renewcommand{\thetlemma}{}
\begin{tlemma} $\rho:\stfs\rightarrow\stfh$ is a Lie superalgebra isomorphism.
\end{tlemma}
{\bf Proof: }Let
$T_{ij}^{\sharp}(a,b)=[X_{ij}^{\sharp}(a),X_{ji}^{\sharp}(b)]$.
Then one can easily check that for $a,b\in R$ and distinct
$i,j,k$, one has
\begin{align}
&T_{ij}^{\sharp}(a,b)=-(-1)^{\omega(i)+\omega(j)}T_{ji}^{\sharp}(b,a)                       \tag{3.15}\\
&T_{ij}^{\sharp}(ab,c)=T_{ik}^{\sharp}(a,bc)+(-1)^{\omega(i)+\omega(k)}T_{kj}^{\sharp}(b,ca).\tag{3.16}
\end{align}
Indeed, the proof of (3.16) is the same as the proof in [MP, Lemma
4.1]. Put
$t^{\sharp}(a,b)=T_{1j}^{\sharp}(a,b)-T_{1j}^{\sharp}(1,ab)$ for
$a,b\in R, 2\leq j\leq 4$. Then $t^\sharp(a,b)$ does not depend on
the choice of $j$. Also, one can easily check (as in the proof of
Lemma 2.15 in [GS]) that
$$\stfs={\frak T}^{\sharp}\oplus_{1\leq i \neq j\leq 4}
X_{ij}^{\sharp}(R)$$ where
$${\frak T}^{\sharp}= \left(\sum_{i, j, k, l \text{ are distinct}}[ X_{ij}^{\sharp}(R), X_{kl}^{\sharp}(R)]\right)
\oplus\left(\sum_{1\leq i < j\leq 4}[ X_{ij}^{\sharp}(R),
X_{ji}^{\sharp}(R)]\right).$$ It then follows from (3.15) and
(3.16) above that
\begin{equation}{\frak T}^{\sharp}={\mathcal W}\oplus\left(t^{\sharp}(R,R)\oplus T_{12}^{\sharp}(1,R)\oplus T_{13}^{\sharp}(1,R)\oplus
T_{14}^{\sharp}(1,R)\right)\tag{3.17}
\end{equation} where
$t^{\sharp}(R,R)$ is the linear span of the elements
$t^{\sharp}(a,b)$. So by Lemma 1.19, it suffices to show that the
restriction of $\rho$ to $t^{\sharp}(R,R)$ is injective.

Now the similar argument as  given in [AG, Lemma 6.18] shows that
 there exists a linear map from $t(R,R)$ to $t^{\sharp}(R,R)$ so
that $t(a,b)\mapsto t^{\sharp}(a,b)$ for $a,b\in R$. This map is
the inverse of the restriction of $\rho$ to $t^{\sharp}(R,R)$.
$\Box$

The following theorem is the main result of this section:
\newtheorem{thmf}{Theorem 3.18}
\renewcommand{\thethmf}{}
\begin{thmf}$(\stfh,\pi)$ is the universal central extension of $\stf$ and
hence
$$H_2(\stf)\cong(0)_{\bar{0}}\oplus({\mathcal
W})_{\bar{1}}.$$
\end{thmf}
{\bf Proof: } We imitate the method of proving the universal central
extension of ${\frak st}_4(R)$ in [GS].

 Suppose that
\begin{equation}0\rightarrow{\mathcal
V}\rightarrow\stft\overset{\tau}{\rightarrow}\stf\rightarrow\notag
0\end{equation} is a central extension of $\stf$. We must show that
there exists a Lie superalgebra homomorphism
$\eta:\stfh\rightarrow\stft$ so that $\tau\circ\eta=\pi$. Thus, by
Lemma 3.14, it suffices to show that there exists a Lie superalgebra
homomorphism $\xi:\stfs\rightarrow\stft$ such that
$\tau\circ\xi=\pi\circ\rho$.

Using the $K$-basis $\{r_\lambda\}_{\lambda\in\Lambda}$ of $R$, we
choose a preimage $\widetilde{X}_{ij}(a)$ of $X_{ij}(a)$ under
$\tau$, $1\leq i\neq j\leq 4,a\in
\{r_\lambda\}_{\lambda\in\Lambda}$, so that the elements
$\widetilde{X}_{ij}(a)$ satisfy the relations (3.6)-(3.12). For
distinct $i ,j, k$, let
$$[\widetilde{X}_{ik}(a),\widetilde{X}_{kj}(b)]=\widetilde{X}_{ij}(ab)+{\mu}_{ij}^k(a,b)$$
where ${\mu}_{ij}^k(a,b)\in{\mathcal V}$. Take distinct $
i,j,k,l$, then
$$\left [\widetilde{X}_{ik}(a),[\widetilde{X}_{kl}(c),\widetilde{X}_{lj}(b)]\right
]=[\widetilde{X}_{ik}(a),\widetilde{X}_{kj}(cb)].$$ But the left
side is, by super-Jacobi identity,
$$\left [[\widetilde{X}_{ik}(a),\widetilde{X}_{kl}(c)],\widetilde{X}_{lj}(b)]\right
]+(-1)^{(\omega(i)+\omega(k))(\omega(k)+\omega(l))}\left
[\widetilde{X}_{kl}(c),[\widetilde{X}_{ik}(a),\widetilde{X}_{lj}(b)]\right
]=[\widetilde{X}_{il}(ac),\widetilde{X}_{lj}(b)].$$ as
$[\widetilde{X}_{ik}(a),\widetilde{X}_{lj}(b)]\in \mathcal V$.
Thus
\begin{equation}
[\widetilde{X}_{ik}(a),\widetilde{X}_{kj}(cb)]=[\widetilde{X}_{il}(ac),\widetilde{X}_{lj}(b)].\notag
\end{equation}In particular, ${\mu}_{ij}^k(a, cb)={\mu}_{ij}^l(ac, b)$ and
$[\widetilde{X}_{ik}(a),\widetilde{X}_{kj}(b)]=[\widetilde{X}_{il}(a),\widetilde{X}_{lj}(b)]$.
It follows that ${\mu}_{ij}^k(a,b)={\mu}_{ij}^l(a,b)={\mu}_{ij}(a,
b)$ which show ${\mu}_{ij}^k(a,b)$ is independent of the choice of
$k$ and ${\mu}_{ij}(c, b)={\mu}_{ij}(1, cb)$, we have
\begin{equation}
[\widetilde{X}_{ik}(a),\widetilde{X}_{kj}(b)]=\widetilde{X}_{ij}(ab)+{\mu}_{ij}(a,b).\notag
\end{equation}
Taking $a=1$, we have
\begin{equation}
[\widetilde{X}_{ik}(1),\widetilde{X}_{kj}(b)]=\widetilde{X}_{ij}(b)+{\mu}_{ij}(1,b).\notag
\end{equation} Now, we replace $\widetilde{X}_{ij}(b)$ by
$\widetilde{X}_{ij}(b)+{\mu}_{ij}(1,b)$. Then the elements
$\widetilde{X}_{ij}(b)$ still satisfy the relations (3.6).
Moreover we have
\begin{equation}
[\widetilde{X}_{ik}(a),\widetilde{X}_{kj}(b)]=\widetilde{X}_{ij}(ab)\tag{3.19}
\end{equation}
for $a,b\in R$ and distinct $i,j,k$. So the elements
$\widetilde{X}_{ij}(a)$ satisfy (3.7).

Next for $k\neq i , k\neq j$, we have
\begin{align}
&[\widetilde{X}_{ij}(a),\widetilde{X}_{ij}(b)]\\
=&\left[\widetilde{X}_{ij}(a),[\widetilde{X}_{ik}(b),\widetilde{X}_{kj}(1)]\right]\notag\\
=&\left[[\widetilde{X}_{ij}(a),\widetilde{X}_{ik}(b)],\widetilde{X}_{kj}(1)\right]+(-1)^{(\omega(i)+\omega(j))(\omega(i)+\omega(k))}
\left[\widetilde{X}_{ik}(b),[\widetilde{X}_{ij}(a),\widetilde{X}_{kj}(1)]\right]\notag\\
=&0+0=0\tag{3.20}
\end{align}
as both $[\widetilde{X}_{ij}(a),\widetilde{X}_{ik}(b)]$ and
$[\widetilde{X}_{ij}(a),\widetilde{X}_{kj}(1)]$ are in $\mathcal V$.
Thus, we get the relation (3.9).

For (3.10), taking $l\notin \{i,j,k\}$
\begin{align}
&[\widetilde{X}_{ij}(a),\widetilde{X}_{ik}(b)]\\
=&\left[\widetilde{X}_{ij}(a),[\widetilde{X}_{il}(b),\widetilde{X}_{lk}(1)]\right]\notag\\
=&\left[[\widetilde{X}_{ij}(a),\widetilde{X}_{il}(b)],\widetilde{X}_{kj}(1)\right]+(-1)^{(\omega(i)+\omega(j))(\omega(i)+\omega(l))}
\left[\widetilde{X}_{il}(b),[\widetilde{X}_{ij}(a),\widetilde{X}_{lk}(1)]\right]\notag\\
=&0+0=0\tag{3.21}
\end{align}
with
$[\widetilde{X}_{ij}(a),\widetilde{X}_{il}(b)],[\widetilde{X}_{ij}(a),\widetilde{X}_{lk}(1)]\in
\mathcal V$. Similarly, we have
\begin{equation}
[\widetilde{X}_{ij}(a),\widetilde{X}_{kj}(b)]=0 \tag{3.22}
\end{equation} for distinct $i,j,k$, which is the relation (3.11).

To verify (3.12)  one needs a few more steps. First, set
$\widetilde{T}_{ij}(a,b)=[\widetilde{X}_{ij}(a),\widetilde{X}_{ji}(b)]$.
The following brackets are easily checked by the super-Jacobi
identity.
\begin{align}[\widetilde{T}_{ij}(a,b),\widetilde{X}_{ik}(c)]=\widetilde{X}_{ik}(abc), & \ \
[\widetilde{T}_{ij}(a,b),\widetilde{X}_{kj}(c)]
=(-1)^{\omega(i)+\omega(j)}\widetilde{X}_{kj}(cba)\notag\\
\text{ and }[\widetilde{T}_{ij}(a,b),\widetilde{X}_{kl}(c)]&=0.
\tag{3.23}
\end{align}
Then we have
\begin{align}
&[\widetilde{T}_{ij}(a,b),\widetilde{X}_{ij}(c)]\\
=&\left[\widetilde{T}_{ij}(a,b),[\widetilde{X}_{ik}(c),\widetilde{X}_{kj}(1)]\right]\notag\\
=&\left[[\widetilde{T}_{ij}(a,b),\widetilde{X}_{ik}(c)],\widetilde{X}_{kj}(1)\right]+
\left[\widetilde{X}_{ik}(c),[\widetilde{T}_{ij}(a,b),\widetilde{X}_{kj}(1)]\right]\notag\\
=&\widetilde{X}_{ij}(abc)+(-1)^{\omega(i)+\omega(j)}\widetilde{X}_{ij}(cba)\\
=&\widetilde{X}_{ij}(abc+(-1)^{\omega(i)+\omega(j)}cba)\tag{3.24}
\end{align}
for $a,b,c\in R$ and distinct $i,j,k,l$.

Next for distinct $i,j,k,l$, let
$$[\widetilde{X}_{ij}(a),\widetilde{X}_{kl}(b)]=\nu^{ij}_{kl}(a,b)$$
where $\nu^{ij}_{kl}(a,b)\in \mathcal V$.

Since one and only one between $\widetilde{X}_{ij}(a)$ and
$\widetilde{X}_{kl}(b)$ is even, we can assume
deg$(\widetilde{X}_{ij}(a))=\bar{0}$. By (3.23) and (3.24),
\begin{align}
2\nu^{ij}_{kl}(a,b)&=[\widetilde{X}_{ij}(2a),\widetilde{X}_{kl}(b)]
=\left[[\widetilde{T}_{ij}(1,1),\widetilde{X}_{ij}(a)],\widetilde{X}_{kl}(b)]\right]\notag\\
&=\left[\widetilde{T}_{ij}(1,1),[\widetilde{X}_{ij}(a),\widetilde{X}_{kl}(b)]\right]+
\left[[\widetilde{T}_{ij}(1,1),\widetilde{X}_{kl}(b)],\widetilde{X}_{ij}(a)]\right]\notag\\
&=0\notag
\end{align}
which yields
\begin{equation}
\nu^{ij}_{kl}(a,b)=-\nu^{ij}_{kl}(a,b).\tag{3.25}
\end{equation} for any distinct $1\leq i,j,k,l\leq 4$ and $a,b\in
R$.

Thus, with the universal property of $\stfs$, we can obtain the Lie
superalgebra homomorphism $\xi:\stfs\rightarrow\stft$ so that
$\tau\circ\xi=\pi\circ\rho$ (as was done in the proof of [GS,
Theorem 2.19]). $\Box$

\noindent {\bf Remark 3.26 }If $2$ is an invertible element of
$K$, then $R=2R$. Thus ${\mathcal I}_2=R$ and ${\mathcal
W}=R_2^6=0$. In this case, $\stf$ is centrally closed.

If the characteristic of $K$ is $2$, we display the following two
examples which are two extreme cases.

\noindent{\bf Example 3.27}  Let $R$ be an associative commutative
$K$-algebra where char $K=2$,  then we have ${\mathcal I}_2 =0$
and $R_2 = R$. Therefore $H_2(\stf)=R^6$.

\noindent{\bf Example 3.28} Let $K$ be a field of characteristic
two. $R=W_k$ is the Weyl algebra which is a unital associative
algebra over $K$ generated by $x_1, \dots, x_k, y_1, \dots, y_k$
subject to the relations $x_ix_j=x_jx_i, \ y_iy_j=y_jy_i, \
x_iy_j-y_jx_i=\delta_{ij}$. Then ${\mathcal I}_2=R$, $H_2(\stf)=0$
and $\stf$ is centrally closed.

\vspace{4mm}

\noindent{\bf \S 4 Central extensions of $\stg$}

In this section, we compute the universal central extension $\stgh$
of $\stg$.

\newtheorem{sdef}{Definition 4.1}
\renewcommand{\thesdef}{}
\begin{sdef}
${\mathcal U}=R_2^4\oplus R_0^2$ is the direct sum of four copies
of $R_2$ and two copies of $R_0$,
$\epsilon_{m}(\bar{a})=(0,\cdots,\overline{a},\cdots,0)$ is the
element of $\mathcal U$, of which the $m$-th component is
$\overline{a}$ and others are zero, for $1\leq m\leq 6$
\end{sdef}

Recall the set of all the quadruple with the distinct components
$P$ and the action of $S_4$ on $P$. $H=\{(1),(13),(24),(13)(24)\}$
is the subgroup of $S_4$, and $P$ have the partition
$P=\bigsqcup_{m=1}^6P_m$, where $P_m=(\sigma_mH)((1,2,3,4))$(cf.
Section 3).

Since the index set of $\stg$ is $\Omega=\{1,2\}\biguplus\{3,4\}$,
we need to classify the subsets $P_m$ of $P$ in this partition.

\newtheorem{sprop}{Proposition 4.2}
\renewcommand{\thesprop}{}
\begin{sprop}
If an element $(i,j,k,l)\in P_m$ satisfies $\omega(i)=\omega(k)$,
then all the elements of $P_m$ have this property.
\end{sprop}
{\bf Proof:} In fact, $\omega(i)=\omega(k)$ induces
$\omega(j)=\omega(l)$. It is easy to see that it is preserved
under the action of $H$. The result is obvious. $\Box$

One can easily see that the ones of $P_m$ satisfying the above
property are:
$$\{(1,3,2,4),(1,4,2,3),(2,3,1,4),(2,4,1,3)\}$$
and
$$\{(3,1,4,2),(3,2,4,1),(4,1,3,2),(4,2,3,1)\}.$$

We denote them by  $P_5$ and $P_6$ respectively. Fix the index map
$\theta:P\rightarrow\{1,2,3,4,5,6\}$ satisfying
$\theta((1,3,2,4))=5$, $\theta((3,1,4,2))=6$. For $(i,j,k,l)\in
P_5\bigsqcup P_6$, we always have $\omega(i)+\omega(j)=\bar{1}$ and
$\omega(k)+\omega(l)=\bar{1}$.

Here we single out $P_5$ and $P_6$ from the others. This is because
for $1\leq m\leq4$, there exist elements $(i,j,k,l)$ of $P_m$ such
that $\omega(i)+\omega(j)=\bar{0}$, but it is not true for $P_5$ and
$P_6$.

Now we take a $K$-basis $\Gamma$ of $\stf$, which contains
$\{X_{ij}(r)|r\in\{r_\lambda\}_{\lambda\in\Lambda},1\leq i\neq
j\leq4\}$. Define $\psi: \Gamma \times \Gamma \to \mathcal W$ by
$$\psi(X_{ij}(r),X_{kl}(s))=sign((i,j,k,l))\epsilon_{\theta((i,j,k,l))}(\overline{rs})$$
 for $r,s\in\{r_\lambda\}_{\lambda\in\Lambda}$, $(i,j,k,l)\in P$ and $\psi=0$, otherwise.

We take the the symbols $sign((i,j,k,l))=1$ for
$(i,j,k,l)\in\bigsqcup_{m=1}^4P_m$, and
$$sign((i,j,k,l))=\cases 1\  &{\text if }\
(i,j,k,l)=(1,3,2,4),(2,4,1,3),(3,1,4,2),(4,2,3,1)\\
-1\ &{\text if }\
(i,j,k,l)=(1,4,2,3),(2,3,1,4),(3,2,4,1),(4,1,3,2)\endcases$$ on
$P_5\bigsqcup P_6$

We then have
\newtheorem{alemma}{Lemma 4.3}
\renewcommand{\thealemma}{}
\begin{alemma} The bilinear map $\psi$ is a (super) $2$-cocycle.
\end{alemma}

\noindent{\bf Proof: }By the definition, one can check the (super)
skew-symmetry of $\psi$.

In fact, if $1\leq m\leq 4$, $\epsilon_m(\bar{a})=
-\epsilon_m(\bar{a})$, thus the $\pm$ sign don't play any role for
$(i,j,k,l)\in\bigsqcup_{m=1}^4P_m$. On the other hand,
$$\psi(X_{ij}(r),X_{kl}(s))=\psi(X_{kl}(r),X_{ij}(s))=-(-1)^{(\omega(i)+\omega(j))(\omega(k)+\omega(l))}\psi(X_{kl}(r),X_{ij}(s)),$$
for $(i,j,k,l)\in P_5\bigsqcup P_6$,
$\omega(i)+\omega(j)=\omega(k)+\omega(l)=\bar{1}$ and the
definition of $sign((i,j,k,l))$.

Moreover, $\psi$ is skew-symmetric on $\stg_{\bar{0}}$ and it is
clear that $\psi(\gamma,\gamma)=0$ for $\gamma$ is contained in the
fixed $K$-basis $\Gamma$ of $\stg_{\bar{0}}$, which implies
$\psi(w,w)=0$, for any $w\in\stg_{\bar{0}}$.

Next, we should show $J(x,y,z)=0$, where
$$J(x,y,z)=(-1)^{{\text deg}(x){\text
deg}(z)}\psi([x,y], z) + (-1)^{{\text deg}(x){\text
deg}(y)}\psi([y, z], x) + (-1)^{{\text deg}(y){\text
deg}(z)}\psi([z, x], y)$$ for the homogenous elements $x,y,z$.
According to Lemma 1.17 and Lemma 1.19, the Steinberg Lie
superalgebra $\stg$ has the decomposition :
\begin{align} \stg=& t(R,R)\oplus T_{12}(1,R)\oplus T_{13}(1,R)\oplus
T_{14}(1,R)\notag \\
&\oplus_{1\leq i \neq j\leq n}X_{ij}(R), \tag{4.4}
\end{align}
where $t(R,R)$ is the $K$-linear span of the elements $t(a,b)$.

We will show the following two possibilities:

{\bf Case 1:} Clearly, the number of elements of ${x,y,z}$
belonging to the subalgebra $\frak T$ such that $\psi([x,y],z)\neq
0$ is at most one. Thus we can suppose that
$x=X_{ij}(a),y=X_{kl}(b)$ and $z\in \frak T$. If
$(i,j,k,l)\in\bigsqcup_{m=1}^4P_m$, it is similar with the proof
of [GS, Lemma2.3]. Therefore, we only should consider
$(i,j,k,l)\in P_5\bigsqcup P_6$. Fix $x=X_{13}(a),y=X_{24}(b)$ and
omit the other similar cases. By (4.4), we can assume that either
$z=t(c,d)$, where $c,d\in R$, or $z=T_{1j}(1,c)$, where $2\leq
j\leq 4$ and $c\in R$. Note that deg$(z)=\bar{0}$,
$\theta((1,3,2,4))=\theta((2,4,1,3))=5$ and
$sign((1,3,2,4))=sign((2,4,1,3))=1$. By Lemma 1.9, when
$z=t(a,b)$, we have
\begin{eqnarray}
J(x,y,z))&=&\psi([t(c,d),X_{13}(a)],X_{24}(b))\nonumber\\
&=&\psi(X_{13}((cd-dc)a),X_{24}(b))\nonumber\\
&=&\epsilon_5(\overline{(cd-dc)ab)}=0;\nonumber
\end{eqnarray}
when $z=T_{12}(c)$,
\begin{eqnarray}
J(x,y,z)&=&-\psi([X_{24}(b),T_{12}(1,c)],X_{13}(a))+\psi([T_{12}(1,c),X_{13}(a)],X_{24}(b))\nonumber\\
&=&\psi([T_{12}(1,c),X_{24}(b)],X_{13}(a))+\psi(X_{13}(ca),X_{24}(b))\nonumber\\
&=&-\psi(X_{24}(cb),X_{13}(a))+\psi(X_{13}(ca),X_{24}(b))\nonumber\\
&=&-\epsilon_5(\overline{cba})+\epsilon_5(\overline{cab})=\epsilon_5(\overline{c(ab-ba)})=0;\nonumber
\end{eqnarray}
when $z=T_{13}(c)$,
\begin{eqnarray}
J(x,y,z)&=&-\psi([X_{24}(b),T_{13}(c)],X_{13}(a))+\psi([T_{13}(1,c),X_{13}(a)],X_{24}(b))\nonumber\\
&=&0+\psi(X_{13}(ac+(-1)^{\omega(1)+\omega(3)}ca),X_{3,4}(b))\nonumber\\
&=&\epsilon_5((\overline{ac-ca)b})=0;\nonumber
\end{eqnarray}
when $z=T_{14}(c)$,
\begin{eqnarray}
J(x,y,z)&=&-\psi([X_{24}(b),T_{14}(c)],X_{13}(a))+\psi([T_{14}(1,c),X_{13}(a)],X_{2,4}(b))\nonumber\\
&=&\psi([T_{14}(c),X_{24}(b)],X_{13}(a))+\psi(X_{13}(ca),X_{24}(b))\nonumber\\
&=&-\psi(X_{24}(bc),X_{13}(a))+\psi(X_{13}(ca),X_{24}(b))\nonumber\\
&=&-\epsilon_5(\overline{bca})+\epsilon_5(\overline{cab})=\epsilon_5(\overline{cab-bca})=0.\nonumber
\end{eqnarray}
{\bf Case 2:} If there is none of $\{x,y,z\}$ belonging to $\frak
T$, the nonzero terms of $J(x,y,z)$ must be
$\psi([X_{ik}(a),X_{kj}(b)],X_{kl}(c))$ or
$\psi([X_{il}(a),X_{lj}(b)],X_{kl}(c))$, for distinct $i,j,k,l$
and $a,b,c\in R$.

If $(i,j,k,l)\in\bigsqcup_{m=1}^4P_m$, it is the same as Case 2 in
the proof of [GS, Lemma 2.3]. Thus, it is enough to check the
following two subcases.

 One is: $x=X_{12}(a)$,$y=X_{23}(b)$,
$z=X_{24}(c)$, and
\begin{eqnarray}
J(x,y,z)&=&\psi(X_{13}(ab),X_{24}(c))-\psi(-X_{14}(ac),X_{23}(b))\nonumber\\
&=&sign((1,3,2,4))\epsilon_{\theta((1,3,2,4))}(\overline{abc})+sign((1,4,2,3))\epsilon_{\theta((1,4,2,3))}(\overline{acb})\nonumber\\
&=&\epsilon_5(\overline{a(bc-cb)})=0.\nonumber
\end{eqnarray}

The other is: $x=X_{14}(a)$,$y=X_{43}(b)$, $z=X_{24}(c)$, and
\begin{eqnarray}
J(x,y,z)&=&-\psi(X_{13}(ab),X_{24}(c))+\psi(-X_{23}(cb),X_{14}(a))\nonumber\\
&=&-sign((1,3,2,4))\epsilon_{\theta((1,3,2,4))}(\overline{abc})-sign((2,3,1,4))\epsilon_{\theta((2,3,1,4))}(\overline{acb}))\nonumber\\
&=&-\epsilon_5(\overline{a(bc-cb)})=0\nonumber
\end{eqnarray}
as $sign((1,3,4,2))=1,sign((1,4,2,3))=sign((2,3,1,4))=-1$ and
$$\theta((1,3,4,2))=\theta((1,4,2,3))=\theta((2,3,1,4))=5$$ for any
$a,b,c\in R$. The proof is completed. $\Box$

\noindent{\bf Remark 4.5}   In view of the proof, for $m=5,6$, the
m-th coordinate doesn't need modular $2R$. In this case, $\psi$ has
already become a (super) $2$-cocycle.

Since $${\mathcal U}=span_{K}\{\psi(X_{ij}(a),X_{kl}(b))|a,b\in R
{\text \ and }\  i,j,k,l {\text \ are \ distinct}\}$$ and ${\text
deg}(X_{ij}(a))={\text deg}(X_{kl}(b))$ for distinct $1\leq
i,j,k,l\leq 4$, we obtain a central extension of Lie superalgebra
$\stg$ satisfying that $\mathcal U$ is the even part of the kernel :

\begin{equation}
0\rightarrow({\mathcal
U})_{\bar{0}}\oplus(0)_{\bar{1}}\rightarrow\stgh\overset{\pi}\rightarrow\stg\rightarrow
0,\tag{4.6}
\end{equation}

i.e. \begin{equation} \stgh=\left(({\mathcal
U})_{\bar{0}}\oplus(0)_{\bar{1}}\right)\oplus\stg.\tag{4.7}
\end{equation}
$(\stgh,\pi)$ is  a central extension of $\stg$. It is similar to
the $\stf$ case, we define a Lie superalgebra $\stgs$ to be the Lie
superalgebra generated by the symbols $X_{ij}^{\sharp}(a)$, $a\in R$
and the $K$-linear space ${\mathcal U}$, with
deg$(X_{ij}^{\sharp}(a))=\omega(i)+\omega(j)$ and deg$(u)=\bar{0}$
for any $u\in{\mathcal U}$, satisfying the following relations:
\begin{align} &a\mapsto X_{ij}{^\sharp}(a) \text{ is a $K$-linear mapping,}\tag{4.8}\\
&[X_{ij}^{\sharp}(a), X_{jk}^{\sharp}(b)] = X_{ik}^{\sharp}(ab), \text{ for distinct } i, j, k, \tag{4.9}\\
&[X_{ij}^{\sharp}(a),{\mathcal U}]=0, \text{ for distinct } i, j, \tag{4.10}\\
&[X_{ij}^{\sharp}(a),X_{ij}^{\sharp}(b)]=0, \text{ for distinct } i, j, \tag{4.11}\\
&[X_{ij}^{\sharp}(a),X_{ik}^{\sharp}(b)]=0, \text{ for distinct } i, j, k, \tag{4.12}\\
&[X_{ij}^{\sharp}(a),X_{kj}^{\sharp}(b)]=0, \text{ for distinct } i, j, k, \tag{4.13}\\
&[X_{ij}^{\sharp}(a),
X_{kl}^{\sharp}(b)]=sign((i,j,k,l))\epsilon_{\theta((i,j,k,l))}(\overline{ab}),
\text{ for distinct } j, k, i, l, \tag{4.14}
\end{align}
where $a,b\in R$,$1\leq i,j,k, l\leq 4$. As $1\in R$, $\stgs$ is
perfect. Clearly, there is a unique Lie superalgebra homomorphism
$\rho:\stgs\rightarrow\stgh$ such that
$\rho(X^\sharp_{ij}(a))=X_{ij}(a)$ and $\rho|_{\mathcal U}=id$.

As was done in Lemma 3.14, we have
\newtheorem{llemma}{Lemma 4.15}
\renewcommand{\thellemma}{}
\begin{llemma} $\rho:\stgs\rightarrow\stgh$ is a Lie superalgebra isomorphism.
\end{llemma}

Now we can state  the main theorem of this section.
\newtheorem{fhmf}{Theorem 4.16}
\renewcommand{\thefhmf}{}
\begin{fhmf}$(\stgh,\pi)$ is the universal central extension of $\stg$ and
hence
$$H_2(\stg)\cong({\mathcal
U})_{\bar{0}}\oplus(0)_{\bar{1}}.$$
\end{fhmf}

\noindent {\bf Proof: }  Suppose that
\begin{equation}0\rightarrow{\mathcal
V}\rightarrow
\stgt\overset{\tau}{\rightarrow}\stg\rightarrow\notag
0\end{equation} is a central extension of $\stg$. We must show
that there exists a Lie algebra homomorphism
$\eta:\stgh\rightarrow\stgt$ so that $\tau\circ\eta=\pi$. Thus, by
Lemma 4.15, it suffices to show that there exists a Lie algebra
homomorphism $\xi:\stgs\rightarrow\stgt$ so that
$\tau\circ\xi=\pi\circ\rho$.

We choose an appropriate  preimage $\widetilde{X}_{ij}(a)$ of
$X_{ij}(a)$ under $\tau$, and check them satisfying (4.8)-(4.14).
The difference from the proof of Theorem 3.18 is to treat
$[\widetilde{X}_{ij}(a),\widetilde{X}_{kl}(b)]$,  which is also
denoted by $\nu^{ij}_{kl}(a,b)$.

We first have
$$\nu^{il}_{kj}(bc,a)=(-1)^{(\omega(k)+\omega(l))(\omega(k)+\omega(j))}\nu^{ij}_{kl}(ba,c).$$ Then taking $b=1$ or $c=1$,
\begin{equation}\nu^{il}_{kj}(b,a)=(-1)^{(\omega(k)+\omega(l))(\omega(k)+\omega(j))}\nu^{ij}_{kl}(a,b)=(-1)^{(\omega(k)+\omega(l))(\omega(k)+\omega(j))}\nu^{ij}_{kl}(ba,1)\tag{4.17}
\end{equation}where $a,b\in R$ and $i,j,k,l$ are distinct.

For $1\leq m\leq 4$, there exists an element $(i,j,k,l)\in P_m$,
such that $\omega(i)+\omega(j)=\bar{0}$, by (3.24), we obtain
\begin{equation}2\nu^{ij}_{kl}(a,b)=0\tag{4.18}
\end{equation}where $a,b\in R$.
As in the proof of Theorem 3.18, one has
\begin{equation}\nu^{ij}_{kl}({\mathcal I}_2,1)=0.\tag{4.19}\end{equation}
By (4.16), the equation holds for any $(i,j,k,l)\in
\bigsqcup_{m=1}^4P_m$.

On the other hand, if $m=5,6$, for all $(i,j,k,l)\in P_m$,
$\omega(i)+\omega(j)=\bar{1}$, then
\begin{align}\nu^{ij}_{kl}(c(ab-ba),1)&=\nu^{ij}_{kl}(ab-ba,c)=\nu^{ij}_{kl}(ab+(-1)^{\omega(i)+\omega(j)}ba,c)\notag\\
&=[\widetilde{X}_{ij}(ab+(-1)^{\omega(i)+\omega(j)}ba),\widetilde{X}_{kl}(c)]\notag\\
&=\left[[\widetilde{T}_{ij}(a,b),\widetilde{X}_{ij}(1)],\widetilde{X}_{kl}(c)\right]\notag\\
&=0\notag
\end{align}for $a,b,c,\in R$, which shows
\begin{equation}\nu^{ij}_{kl}({\mathcal I}_0,1)=0\tag{4.20}\end{equation}
for $(i,j,k,l)\in P_5\bigsqcup P_6$.

The rest of the proof is similar to Theorem 3.18, we can obtain
$\xi:\stgs\rightarrow\stgt$. The only difference is that we need
paying attention to the sign of the restriction of $\xi$ on
$\mathcal U$ as the $5$-th and $6$-th coordinate component of
$\mathcal U$ is $R_0$. Let
$\xi(\epsilon_m(\bar{a}))=sign((i,j,k,l))\nu^{ij}_{kl}(1,a)$, where
$sign(i,j,k,l)$ is defined before Lemma 4.3. It is easy to see that
the choice of sign coincides with the (super) skew-symmetry and
(4.17). Thus, the Lie homomorphism $\psi$ is well defined on
$\mathcal U$.
 $\Box$

\noindent {\bf Remark 4.20 }Note that $H_2(\stg)\cong R_2^4\oplus
R_0^2$. Even $2$ is an invertible element of $K$ so that $R=2R$ and
$R_2=0$, $R_0$ is not necessarily equal to $0$. Particularly, if $R$
is commutative, then ${\mathcal I}_0=R[RR]=0$ and $R_0=R$. In this
case, $H_2(\stg)\cong R^2$ which is not trivial.

\vspace{4mm}

\noindent{\bf \S 5 Concluding remarks}

Combining Theorem 1.19, Theorem 2.1, Theorem 3.18 and Theorem
4.15, we completely determined $H_2(\st)$ for $m+n\geq 3$.

\newtheorem{ghmf}{Theorem 5.1}
\renewcommand{\theghmf}{}
\begin{ghmf} let $K$ be a unital commutative ring and $R$ be a
unital associative $K$-algebra. Assume that $R$ has a $K$-basis
containing the identity element. Then
\begin{equation}
H_2(\st)=\cases 0 &\text{ for } m+n=3 \ {\text and } \ m+n\geq 5\\
(0)_{\bar{0}}\oplus (R_2^6)_{\bar{1}}  &\text{ for } m=3,n=1\\
(R_2^4\oplus R_0^2)_{\bar{0}}\oplus(0)_{\bar{1}} &\text{ for }
m=2,n=2\endcases\notag
\end{equation}
which are ${\mathbb Z}_2$-graded spaces.

\end{ghmf}

It then follows from [MP] that

\newtheorem{ihmf}{Theorem 5.2}
\renewcommand{\theihmf}{}
\begin{ihmf} let $K$ be a unital commutative ring and $R$ be a
unital associative $K$-algebra. Assume that $R$ has a $K$-basis
containing the identity element. Then
\begin{equation}
H_2(sl_n(R))=\cases (HC_1(R))_{\bar{0}}\oplus(0)_{\bar{1}} &\text{ for } m+n=3 \ {\text and } \ m+n\geq 5\\
 (HC_1(R))_{\bar{0}}\oplus (R_2^6)_{\bar{1}} &\text{ for } m=3,n=1\\
(R_2^4\oplus R_0^2\oplus HC_1(R))_{\bar{0}}\oplus(0)_{\bar{1}}
&\text{ for } m=2,n=2\endcases\notag
\end{equation}
where $HC_1(R)$ is the first cyclic homology group of the
associative $K$-algebra $R$ (See [L]).

\end{ihmf}

Department of Mathematics

University of Science and Technology of China

Hefei, Anhui

P. R. China  230026

hjchen@@mail.ustc.edu.cn,

\

Department of Mathematics and Statistics

York University

Toronto, Ontario

Canada  M3J 1P3

ygao@@yorku.ca

and

Department of Mathematics

University of Science and Technology of China

Hefei, Anhui

P. R. China  230026

skshang@@mail.ustc.edu.cn


\begin{thebibliography}{ABG}

\bibitem[ABG]{ }  B. N. Allison, G. M. Benkart, Y. Gao, {\it Central extensions of Lie algebras graded by finite
root systems,} Math. Ann.  316 (2000) 499--527.

\bibitem[AF]{  } B. N. Allison and J. R. Faulkner, {\it Nonassociative coefficient algebras for Steinberg unitary Lie algebras,}
 J. Algebra 161 (1993) 1--19.

\bibitem[AG]{  } B. N. Allison and Y. Gao, {\it Central quotients and corverings of Steinberg unitary Algebras,}
 Canad. J. Math. 17 (1996),261--304.

\bibitem[BeM]{ } G. M. Benkart and R. V. Moody, {\it Derivations, central extensions and affine Lie algebras,}
 Algebras, Groups and Geometries 3 (1986) 456--492.


\bibitem[Bl]{  } S. Bloch, {\it The dilogarithm and extensions of Lie algebras,} Alg. K-theory, Evanston 1980,
Springer Lecture Notes in Math 854 (1981) 1--23.

\bibitem[G1]{  } Y. Gao, {\it Steiberg Unitary Lie Algebras and Skew-Dihedral Homology,}  J.Algebra,17 (1996),261--304.

\bibitem[G2]{  } Y. Gao, {\it On the Steinberg Lie algebras $st_2(R)$,} Comm. in Alg. 21 (1993) 3691--3706.

\bibitem[GS]{  } Y. Gao and S. Shang, {\it Universal coverings of Steinberg Lie
algebras of small characteristic,} math.QA/0512188.

\bibitem[K]{  } V. Kac, {\it Lie superalgebras,} Adv. Math., 26, No. 1, (1977) 8--96.

\bibitem[Ka]{  } C. Kassel, {\it K\"ahler differentials and coverings of complex simple Lie algebras
extended over a commutative ring,} J. Pure and Appl. Alg. 34
(1984) 265--275.

\bibitem[KL]{  } C. Kassel and J-L. Loday, {\it Extensions centrales d'alg\`ebres de Lie,} Ann. Inst. Fourier 32 (4) (1982) 119--142.

\bibitem[L]{  } J-L. Loday, {\it Cyclic homology,} Grundlehren der mathematischen Wissenschaften 301, Springer 1992.

\bibitem[MP]{  } A. V. Mikhalev and I. A. Pinchuk, {\it Universal central extensions of the matrix Lie superalgebras $sl(m,n,A)$,}
Int. Conf. in H.K.U., AMS, (2000) 111--125.

\bibitem[N]{  } E. Neher, {\it An introduction to universal central extensions of Lie superalgebras,}
 Groups, rings, Lie and Hopf algebras (St. John's, NF, 2001), 141--166 Math. Appl., 555,
 Kluwer Acad. Publ., Dordrecht, 2003.

\end{thebibliography}
\end{document}